\documentclass{amsart}
\usepackage{hyperref}
\usepackage{fullpage}

\input xypic

\newcommand{\CA}        {{\mathcal{A}}}
\newcommand{\CB}        {{\mathcal{B}}}
\newcommand{\CC}        {{\mathcal{C}}}
\newcommand{\CD}        {{\mathcal{D}}}
\newcommand{\CE}        {{\mathcal{E}}}
\newcommand{\CF}        {{\mathcal{F}}}
\newcommand{\CG}        {{\mathcal{G}}}
\newcommand{\CI}        {{\mathcal{I}}}
\newcommand{\CM}        {{\mathcal{M}}}
\newcommand{\CU}        {{\mathcal{U}}}
\newcommand{\CV}        {{\mathcal{V}}}

\newcommand{\Ab}        {\operatorname{Ab}}
\newcommand{\Stab}      {\operatorname{Stab}}
\newcommand{\Ho}        {\operatorname{Ho}}

\newcommand{\ann}       {\operatorname{ann}}
\newcommand{\img}       {\operatorname{image}}
\newcommand{\proj}      {\operatorname{proj}}
\newcommand{\supp}      {\operatorname{supp}}
\newcommand{\zar}       {\operatorname{zar}}

\newcommand{\locid}[1]  {\text{locid}\langle #1\rangle}
\newcommand{\colocid}[1]{\text{colocid}\langle #1\rangle}
\newcommand{\loc}[1]    {\operatorname{loc}\langle #1\rangle}
\newcommand{\thickid}[1]{\operatorname{thickid}\langle #1\rangle}

\newcommand{\colim}  {\operatornamewithlimits{\underset{\longrightarrow}{lim}}}
\newcommand{\invlim} {\operatornamewithlimits{\underset{\longleftarrow}{lim}}}
\newcommand{\pgldim}    {\operatorname{pgldim}}

\newcommand{\N}         {{\mathbb{N}}}
\newcommand{\Q}         {{\mathbb{Q}}}
\newcommand{\Z}         {{\mathbb{Z}}}
\newcommand{\Fp}        {{\mathbb{F}_p}}

\newcommand{\bId}       {\overline{\text{Id}}}
\newcommand{\ot}        {\otimes}
\newcommand{\st}        {\;|\;}
\newcommand{\opp}       {{\text{op}}}
\newcommand{\tm}        {\times}
\newcommand{\xra}       {\xrightarrow}
\newcommand{\Sg}        {\Sigma}
\newcommand{\sse}       {\subseteq}
\newcommand{\Smash}     {\wedge}
\newcommand{\bigWedge}  {\bigvee}
\newcommand{\kp}        {\kappa}
\newcommand{\Sgi}       {\Sigma^\infty}
\newcommand{\Lm}        {\Lambda}
\newcommand{\hC}        {\widehat{C}}
\newcommand{\hL}        {\widehat{L}}
\newcommand{\op}        {\oplus}
\newcommand{\xla}       {\xleftarrow}
\newcommand{\Wedge}     {\vee}
\newcommand{\bc}[1]     {\langle #1\rangle}
\newcommand{\pri}       {\mathfrak{p}}

\renewcommand{\:}{\colon}

\newtheorem{theorem}{Theorem}[section]

\theoremstyle{definition}

\newtheorem{definition}[theorem]{Definition}

\numberwithin{equation}{subsection}

\begin{document}
\title{Axiomatic stable homotopy --- a survey}
\author{N.~P.~Strickland}
\address{
Department of Pure Mathematics\\
University of Sheffield\\
Sheffield S3 7RH\\
UK
}
\email{N.P.Strickland@sheffield.ac.uk}
\date{\today}
\bibliographystyle{habbrv}



\begin{abstract}
 We survey various approaches to axiomatic stable homotopy theory,
 with examples including derived categories, categories of
 (possibly equivariant or localized) spectra, and stable categories of
 modular representations of finite groups.  We focus mainly on
 representability theorems, localisation, Bousfield classes, and
 nilpotence.
\end{abstract}

\maketitle 

\section{Introduction}

Axiomatic stable homotopy theory is the study of triangulated
categories formally similar to the homotopy category of spectra in the
sense of Boardman~\cite{ad:shg,ma:ssa}.  While various authors have
used different systems of axioms, there is broad agreement about the
main examples, of which the following are a sample:
\begin{enumerate}
 \item Boardman's category itself, which we denote by $\CB$.
 \item The subcategories of $E(n)$-local and $K(n)$-local
  spectra~\cite{ra:lrc,host:mkl}.
 \item The homotopy category $\CB_G$ of $G$-spectra (indexed on a
  complete universe), where $G$ is a compact Lie
  group~\cite{lemast:esh,ma:ehc}.
 \item The homotopy category of $A$-modules, where $A$ is a
  commutative ring spectrum.  (Here and throughout this paper, the
  phrase ``ring spectrum'' refers to a strictly associative monoid in
  a suitable geometric category of spectra, such as that defined
  in~\cite{ekmm:rma}.) 
 \item The derived category $\CD_R$ of $R$-modules, for a commutative
  ring $R$.  If we let $HR$ denote the associated Eilenberg-MacLane
  ring spectrum, then $\CD_R\simeq\CD_{HR}$ (by~\cite[Theorem
  IV.2.4]{ekmm:rma}), so this is a special case of the previous
  example.
 \item The stable category $\Stab_{kG}$ of $kG$-modules and projective
  equivalence classes of homomorphisms, where $G$ is a finite group
  and $k$ is a field~\cite{ri:ims}.  This occurs as a quotient of the
  category $\CD_{kG}$, in which the morphism sets are given by group
  cohomology; the morphism sets in $\Stab_{kG}$ itself are more
  closely related to Tate cohomology.
 \item The derived category (in a suitable sense) of
  $MU_*MU$-comodules~\cite{ho:htc}.  
\end{enumerate}

There are some further examples that satisfy some authors' axioms but
not others, or where the axioms have not yet been checked (to the best
of my knowledge):
\begin{enumerate}
 \item The derived category of modules over a noncommutative ring.
 \item The derived category of quasicoherent sheaves over a nonaffine
  scheme.
 \item The homotopy category of $G$-spectra indexed by an incomplete
  $G$-universe. 
 \item Various versions of the category of motivic spectra.  (Motivic
  spaces are discussed in~\cite{move:aht}; at the time of writing,
  there is no published account of the corresponding category of
  spectra.) 
\end{enumerate}
It seems an important problem to decide which of the usual axioms
apply to the motivic stable category, and to see what the axiomatic
literature teaches us about this example.

The main topics that have been discussed from an axiomatic point of
view are as follows.
\begin{itemize}
 \item[(a)] Theorems saying that certain (covariant or contravariant)
  functors are representable, generalising the Brown representability
  theorems for (co)homology theories on (finite or infinite) spectra.
 \item[(b)] Phantom maps.
 \item[(c)] Various kinds of localisation, generalising Bousfield's theory
  of localisation with respect to a homology theory.  Special cases
  such as finite, cofinite, algebraic or smashing localisations.
 \item[(d)] The lattice of Bousfield classes, and various related
  lattices (some of them conjecturally identical to the Bousfield
  lattice).
 \item[(e)] Nilpotence theorems in the spirit of Devinatz, Hopkins and
  Smith.
 \item[(f)] Picard groups, Grothendieck rings, and Euler
  characteristics~\cite{fr:ggs,ma:att,ma:pgg,krre:esh}.
 \item[(g)] Projective classes, and generalisations of the Adams
  spectral sequence~\cite{ch:itc}.
 \item[(h)] Duality theorems generalising those of Verdier and
  Gross-Hopkins~\cite{ne:gdt,fahuma:ilr}.
\end{itemize}
For~(f) to~(h) we refer the reader to the cited papers and their
bibliographies.  This survey will concentrate on~(a) to~(e).

The relevant literature consists partly of papers that are explicitly
axiomatic, and partly of papers that are nominally restricted to some
particular category, but whose methods allow straightforward
generalisation to other examples.  Some authors are as follows:
\begin{enumerate}
 \item Margolis's book~\cite{ma:ssa} treats $\CB$ from an axiomatic
  point of view; this was an important inspiration for much of the
  later work.  Earlier still, there were relevant papers by Freyd,
  Heller and Joel Cohen.
 \item Neeman has written extensively, particularly on questions
  related to representability and
  localisation~\cite{bone:hlt,chkene:fbr,ne:brd,ne:ctd,ne:dcs,ne:nat,ne:obc,ne:tba,ne:tc}.
  Some papers are restricted to the case of $\CD_R$; often $R$ need
  not be commutative or noetherian.  When working axiomatically, he
  has generally assumed that his triangulated category $\CC$ is
  ``compactly generated'', but not that $\CC$ has a symmetric monoidal
  structure.  The class of categories considered is thus rather large,
  but unfortunately it is not closed under Bousfield localisation.
  Recently he has introduced the more complex notion of a
  ``well-generated'' triangulated category to repair this problem.
 \item Krause has also written extensively on representability,
  localisation, and versions of the Bousfield
  lattice~\cite{bekr:gim,bekr:pis,kr:dts,kr:sst,kr:smc,kr:nwg,kr:brf,kr:brt,kr:cfs,krre:esh,krzw:seg}.
 \item Beligiannis has written a long paper~\cite{be:rha} covering
  many themes in axiomatic stable homotopy, considered as an analog of
  relative homological algebra in the context of triangulated
  categories. 
 \item Benson, Carlson, Rickard and Gnacadja (working in various
  combinations) have proved many results about the categories
  $\Stab_{kG}$ and $\CD_{kG}$, often using methods that transfer
  easily to an axiomatic
  setting~\cite{be:idm,becari:cvi,becari:cviii,becari:tss,begn:pmpi,begn:pmpii,be:pmpiii,bekr:gim,bekr:pis,ri:blr,ri:ims,ri:mtd,ri:ram,ri:sed}.
  Benson and Wheeler have interpreted the Green correspondence in this
  context~\cite{bewh:gci}.
 \item Hovey, Schwede and Shipley have worked in a more rigid context,
  studying Quillen model categories $\CC_0$ such that the homotopy
  category $\Ho(\CC_0)$ is
  triangulated~\cite{ho:mc,sc:shc,scsh:smc,sh:mus}.
 \item May and coauthors have studied the equivariant stable
  categories $\CB_G$, often using methods that transfer easily to an
  axiomatic setting.  This applies particularly to their work on
  duality, traces, and Picard groups~\cite{falema:pge,ma:pgg,ma:att,lemast:esh}.
 \item Hovey and Palmieri and Strickland wrote a
  memoir~\cite{hopast:ash} on axiomatic stable homotopy theory.  We
  assumed much more than Neeman, and thus could obtain results closer
  to those previously known for $\CB$.  In particular, we assume that
  $\CC$ has a closed symmetric monoidal structure.
\end{enumerate}

\section{Axioms}

We next discuss the various axioms that have been used.  We start with
a category $\CC$.

\subsection{Basics}

The category $\CC$ should be triangulated, and should have coproducts
for all families of objects (indexed by a set).  These are core
axioms, used by almost all authors.  Existence of coproducts should be
seen as an important test of the correctness of the technical details
of the definition of $\CC$.  Boardman's category itself came after
several attempts to define a good category of spectra, and it was the
first to be triangulated and coproduct-complete; it rapidly became
clear that Boardman's version was much more convenient than all the
others.  Similarly, the earliest versions of $\CD_R$ incorporated
various boundedness conditions, and so were not coproduct-complete.
Bokstedt and Neeman~\cite{bone:hlt} adjusted the definitions to remove
this problem, and this allowed much smoother comparisons between
$\CD_R$ and $\CB$.

We recall the definition of a triangulation:
\begin{definition}\label{defn-triangulated-category}
 A \emph{triangulation} of an additive category $\CC$ is an additive
 (\emph{suspension}) functor $\Sigma\:\CC\xra{}\CC$ giving an
 automorphism of $\CC$, together with a collection
 $\boldsymbol{\triangle}$ of diagrams, called \emph{distinguished
   triangles} or \emph{cofibre sequences}, of the form
 \[ X \xra{} Y \xra{} Z \xra{}\Sigma X\]
 such that
 \begin{itemize}
 \item [1.] Any diagram isomorphic to a cofibre sequence is a cofibre
  sequence.
 \item [2.] Any diagram of the following form is a cofibre sequence:
  \[ 0 \xra{} X \xra{1} X \xra{}0 \]
 \item [3.] The first of the following diagrams a cofibre sequence iff
  the second is a cofibre sequence:
  \[ X \xra{f}Y\xra{g}Z\xra{h}\Sigma X\]
  \[ Y\xra{g}Z\xra{h}\Sigma X \xra{-\Sigma f} \Sigma Y. \]
 \item [4.] For any map $f\:X\xra{} Y$, there is a cofibre sequence of
  the following form:
  \[ X\xra{f}Y\xra{}Z\xra{}\Sigma X\]
 \item [5.] Suppose we have a diagram as shown below (with $h$
  missing), in which the rows are cofibre sequences and the
  rectangles commute.  Then there exists a (nonunique) map $h$ making
  the whole diagram commutative.
  \[
   \setlength{\unitlength}{0.07em}
   \begin{picture}(300,80)
   \put(  0,  0){\makebox(0,0){$X$}}
   \put( 80,  0){\makebox(0,0){$Y$}}
   \put(160,  0){\makebox(0,0){$Z$}}
   \put(240,  0){\makebox(0,0){$\Sigma X$}}
   \put(  0, 80){\makebox(0,0){$U$}}
   \put( 80, 80){\makebox(0,0){$V$}}
   \put(160, 80){\makebox(0,0){$W$}}
   \put(240, 80){\makebox(0,0){$\Sigma U$}}
   \put( 11,  0){\vector( 1, 0){58}}
   \put( 87,  0){\vector( 1, 0){58}}
   \put(167,  0){\vector( 1, 0){55}}
   \put(  7, 80){\vector( 1, 0){58}}
   \put( 87, 80){\vector( 1, 0){58}}
   \put(167, 80){\vector( 1, 0){55}}
   \put(  0, 68){\vector( 0,-1){56}}
   \put( 80, 68){\vector( 0,-1){56}}
   \put(240, 68){\vector( 0,-1){56}}
   \multiput(160, 68)( 0,-10){6}{\line( 0,-1){5}}
   \put(160, 17){\vector( 0,-1){ 5}}
   \put(167, 40){\makebox(0,0){$h$}}
   \put(7, 40){\makebox(0,0){$f$}}
   \put(252, 40){\makebox(0,0){$\Sigma f$}}
  \end{picture} \]
 \item [6.] Verdier's octahedral axiom holds: Suppose we have maps
  $X\xra{v}Y\xra{u}Z$, and cofibre sequences $(X,Y,U)$,
  $(X,Z,V)$ and $(Y,Z,W)$ as shown in the diagram.  (A circled arrow
  $U\longrightarrow\hspace{-1.3em}\circ\hspace{0.8em} X$ means a map
  $U\xra{}\Sigma X$.)  Then there exist maps $r$ and $s$ as shown,
  making $(U,V,W)$ into a cofibre sequence, such that the
  following commutativities hold:
  \[ au = rd \qquad \qquad es = (\Sigma v) b \qquad\qquad
     sa = f  \qquad \qquad br = c \]
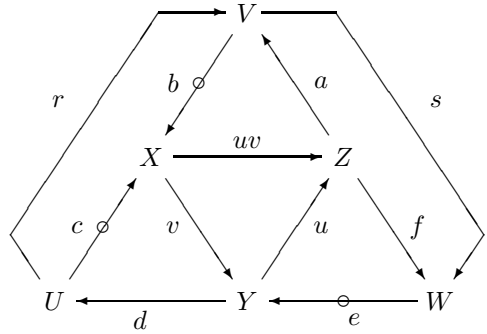
  \[
  \setlength{\unitlength}{0.14em}
  \begin{picture}(140,110)(0,-10)
   \put( 18,11){\vector( 2, 3){18}}
   \put( 70,11){\vector( 2, 3){18}}
   \put( 62,77){\vector(-2,-3){18}}
   \put( 44,38){\vector( 2,-3){18}}
   \put( 96,38){\vector( 2,-3){18}}
   \put( 88,50){\vector(-2, 3){18}}
   \put(112, 5){\vector(-1, 0){40}}
   \put( 60, 5){\vector(-1, 0){40}}
   \put( 46,44){\vector( 1, 0){40}}
   \put( 10,11){\line(-2, 3){ 8}}
   \put(  2,23){\line( 2, 3){40}}
   \put( 42,83){\vector( 1, 0){18}}
   \put( 70,83){\line( 1, 0){20}}
   \put( 90,83){\line( 2,-3){40}}
   \put(130,23){\vector(-2,-3){ 8}}
   \put( 27,25){\circle{3}}
   \put( 53,64){\circle{3}}
   \put( 92, 5){\circle{3}}
   \put( 14, 5){\makebox(0,0){$U$}}
   \put( 66, 5){\makebox(0,0){$Y$}}
   \put(118, 5){\makebox(0,0){$W$}}
   \put( 40,44){\makebox(0,0){$X$}}
   \put( 92,44){\makebox(0,0){$Z$}}
   \put( 66,83){\makebox(0,0){$V$}}
   \put( 86,25){\makebox(0,0){$u$}}
   \put( 46,25){\makebox(0,0){$v$}}
   \put( 66,48){\makebox(0,0){$uv$}}
   \put( 86,64){\makebox(0,0){$a$}}
   \put( 46,64){\makebox(0,0){$b$}}
   \put( 20,25){\makebox(0,0){$c$}}
   \put( 37, 0){\makebox(0,0){$d$}}
   \put( 95, 0){\makebox(0,0){$e$}}
   \put(112,25){\makebox(0,0){$f$}}
   \put( 15,59){\makebox(0,0){$r$}}
   \put(117,59){\makebox(0,0){$s$}}
  \end{picture}
  \]
  (If $u$ and $v$ are inclusions of CW spectra, this essentially just
  says that $(Z/X)/(Y/X) = Z/Y$.  The diagram can be turned into an
  octahedron by lifting the outer vertices and drawing an extra line
  from $W$ to $U$.)
 \end{itemize}
\end{definition}

Following the standard topological notation, we write $[X,Y]$ for the
set $\CC(X,Y)$ (of morphisms in $\CC$ from $X$ to $Y$).  We also put
$[X,Y]_n=[\Sg^nX,Y]$. 

If $\CC_0$ is a pointed Quillen model category, then the homotopy
category $\CC=\Ho(\CC_0)$ automatically has structure close to that
described above, except that the functor $\Sg\:\CC\xra{}\CC$ need not
be an equivalence.  Following Hovey, we say that $\CC_0$ is a
\emph{stable model category} if $\Sg$ is an equivalence; if so, one
can show that $\CC$ is triangulated~\cite[Chapter 7]{ho:mc}.  Similar
results appear in~\cite{ma:att}.  Schwede and Shipley have
shown~\cite{scsh:smc} that most stable model categories are
Quillen-equivalent to $\CD_A$ for some ring spectrum $A$, and that the
general case is only a little more general.

Various modifications and refinements of triangulated categories have
been considered by Neeman~\cite{ne:nat}, May~\cite{ma:att},
Franke~\cite{fr:utc} and probably others.  It seems likely that these
all follow from the existence of an underlying model category $\CC_0$
as above.  For the categories $\CC$ occurring in practice, it seems
that there is always an underlying model category, and that any two
natural choices are Quillen equivalent.  It thus seems reasonable to
assume whenever convenient that one is given a category $\CC_0$.

\subsection{Smash products}

In~\cite{hopast:ash}, we assume that our categories $\CC$ come
equipped with a symmetric monoidal product.  We use notation coming
from topology, and thus write $X\Smash Y$ for the monoidal product of
$X$ and $Y$, and $S$ for the unit, so $S\Smash X=X=X\Smash S$.  We
also assume that there are adjoint function objects $F(Y,Z)$, so
$[X,F(Y,Z)]\simeq[X\Smash Y,Z]$, naturally in all variables.  We write
$S^n=\Sg^nS$, so $S^n\Smash S^m=S^{n+m}$ and $\Sg^nX=S^n\Smash X$.
We also put $\pi_nX=[S^n,X]$.

This structure certainly exists in all the categories mentioned so
far, except for the category $\CD_A$ when $A$ is not commutative.
However, the theory of blocks in group algebras decomposes
$\Stab_{kG}$ as a product of smaller categories, which need not have a
symmetric monoidal structure.  There may be similar examples related
to $\CB_G$.

It is natural to require that the smash product be compatible with the
triangulation.  In~\cite{hopast:ash}, we wrote down the most obvious
compatibility conditions:
\begin{enumerate}
 \item The smash product commutes with suspension, so
  $\Sg(X\Smash Y)\simeq X\Smash\Sg Y$.
 \item The functors $X\Smash(-)$ and $F(X,-)$ preserve cofibre
  triangles.  The contravariant functors $F(-,Y)$ preserve cofibre
  triangles up to a sign change.
 \item The twist map $S^1\Smash S^1\xra{}S^1\Smash S^1$ is
  multiplication by $-1$. 
\end{enumerate}
However, May has given strong evidence that further conditions should
be added.  To explain this, consider a map $f\:X\xra{}X$.  Under
suitable finiteness conditions, this has a trace $\tau(f)\:S\xra{}S$.
If we have a cofibre sequence 
\[ X_0\xra{}X_1\xra{}X_2\xra{}\Sg X_0 \]
and compatible maps $f_i\:X_i\xra{}X_i$, it is natural to hope that
$\tau(f_0)-\tau(f_1)+\tau(f_2)=0$; this is suggested by the theory of
Lefschetz numbers, among other things.  It turns out that the
statement must be adjusted slightly: given $f_0$ and $f_1$, one can
choose a compatible $f_2$ such that $\tau(f_0)-\tau(f_1)+\tau(f_2)=0$,
but this may not be the case for all compatible $f_2$'s.  This (and
various extensions) can be proved in Boardman's category, but the
proof cannot be transferred to the axiomatic setting without adding
some more conditions.  In outline, consider two cofibre sequences
\begin{align*}
 X_0 \xra{} X_1 &\xra{} X_2 \xra{} \Sg X_0 \\
 Y_0 \xra{} Y_1 &\xra{} Y_2 \xra{} \Sg Y_0.
\end{align*}
From these we obtain a $4\tm 4$ diagram with vertices $X_i\Smash X_j$,
in which all squares commute, except that one square anticommutes.
By writing in the diagonal composite in each square, we get 18
commutative triangles.  Each commutative triangle fits in an
octahedron, as in Verdier's axiom.  The 18 resulting octahedra have
many vertices and edges in common, and one can hope to add in some
extra vertices and edges making everything fit together more
coherently.  May~\cite{ma:att} has formulated three axioms about this
situation, and explained how they can be checked when $\CC=\Ho(\CC_0)$
for some Quillen model category $\CC_0$.

There are many interesting cases of noncommutative rings (or ring
spectra) $R$ and $R'$ for which $R\not\simeq R'$ but
$\CD_R\simeq\CD_{R'}$; this is a natural extension of Morita theory.
Examples come from Koszul duality, the Fourier-Mukai transform for
sheaves on abelian varieties, tilting complexes in representation
theory, and so on~\cite{ri:mtd,ri:sed,br:etc,dush:ktd,scsh:csm}.
Nonetheless, it seems that there are no examples of this type where
the categories involved have smash products and the equivalence
respects them.  We know of no rigorous results in this direction,
however. 

\subsection{Generation}

A key feature of Boardman's category is that every spectrum $X$ has a
cell structure, and (essentially equivalently) that if $[S^0,X]_*=0$
then $X=0$.  More generally, suppose we have a set $\CG$ of objects in
a triangulated category $\CC$.  We say that $\CG$ \emph{generates}
$\CC$ if there is no proper triangulated subcategory $\CC'\subset\CC$
closed under all coproducts such that $\CG\sse\CC'$.  We also say that
$\CG$ \emph{detects} $\CC$ if every object $X\in\CC$ with $[A,X]_*=0$
for all $A\in\CG$ actually has $X=0$.  (By taking
$\CC'=\{A\st [A,X]_*=0\}$, we see that generation implies detection,
and the converse holds under suitable finiteness conditions.)
Generating sets for the main examples are as follows.
\begin{enumerate}
 \item $\{S^0\}$ generates $\CB$, and $L_ES^0$ generates the
  subcategory of $E$-local spectra~\cite[Section 3.5]{hopast:ash}.
 \item If $X$ is a finite spectrum of type $n$ (in the usual chromatic
  sense) then $L_{K(n)}X$ generates the category of $K(n)$-local
  spectra~\cite[Theorem 7.3]{host:mkl}, and this is generally a better
  choice of generator than $L_{K(n)}S^0$.
 \item The $G$-spectra $G/H_+$ (as $H$ runs over conjugacy classes of
  closed subgroups) generate $\CB_G$.
 \item If $R$ is a ring, then $R$ generates $\CD_R$.  This also works
  for (strictly associative) ring spectra.
 \item The nonprojective simple $kG$-modules generate $\Stab_{kG}$.
\end{enumerate}

All authors in axiomatic stable homotopy theory assume that $\CC$ is
generated by some set $\CG$ of objects, and impose some smallness
conditions on $\CG$.  The details vary between authors, however.

One popular condition is as follows.  We say that an object $A\in\CC$
is \emph{small} (or \emph{compact}) if the natural map
\[ \bigoplus_{i\in I}[A,X_i] \xra{} [A,\bigoplus_i X_i] \]
is an isomorphism for all families of objects $\{X_i\}$.  For example:
\begin{enumerate}
 \item In $\CB$, the small objects are those of the form $\Sg^dX$
  where $d\in\Z$ and $X$ is a finite CW complex.
 \item In the category of $E(n)$-local spectra, the small objects are
  those that can be written as a retract of $\Sg^dL_{E(n)}X$ for some
  $d\in\Z$ and some finite CW complex $X$.
 \item In the category of $K(n)$-local spectra, the small objects are
  those that can be written as a retract of $\Sg^dL_{K(n)}X$ for some
  $d\in\Z$ and some finite CW complex $X$ of type $n$.
 \item In $\CB_G$, the small objects are
  those that can be written as a retract of $\Sg^dX$ for some $d\in\Z$
  and some finite $G$-CW complex $X$.
 \item In $\CD_R$, the small objects are the finite complexes of
  finitely generated projective modules.
 \item In $\Stab_{kG}$, the small objects are the $kG$-modules $M$
  that are finite-dimensional over $k$.
\end{enumerate}
(Most of these facts are proved in~\cite{hopast:ash}, for example.)

We say that $\CC$ is \emph{compactly generated} if there is a set
$\CG$ of small objects that generates $\CC$.  In the terminology
of~\cite{hopast:ash}, a stable homotopy category is \emph{algebraic}
iff it is compactly generated.  This is a very convenient condition,
and is often satisfied; in particular, the generators listed above for
$\CB$, $\CB_G$, $\CD_R$ and $\Stab_{kG}$ are all small.  In the case
of the $K(n)$-local category, the obvious generator $L_{K(n)}S^0$ is
\emph{not} small, but $L_{K(n)}X$ is small whenever $X$ is finite of
type $n$, so the category is nonetheless compactly generated.  For a
simpler example of the same phenomenon, let $\CC$ be the
$p$-completion of the category $\CD_{\Z}$; then the obvious generator
is $\Z_p$ (which is not small), but the object $\Z/p$ is also a
generator, and is small.  (There is a well-understood axiomatic
framework covering both of these examples: see~\cite[Section
3.3]{hopast:ash}, and Section~\ref{sec-special} of the present paper.)
More seriously, there are many known spectra $E$ for which the
category $\CC_E$ of $E$-local spectra has no nontrivial small objects,
so in particular, $\CC_E$ is not compactly
generated~\cite{ne:ncg}\cite[Appendix B]{host:mkl}.  For example, this
applies with $E=BP$ or $E=H$ or $E=\bigWedge_{0\leq n<\infty}K(n)$.

Another useful condition is \emph{dualisability}.  To formulate this,
we need to assume that $\CC$ has a symmetric monoidal smash product,
as in the previous section.  We write $DA=F(A,S)$, where $S$ is the
unit object for the smash product.  We say that $A$ is
\emph{dualisable} if the natural map $DA\Smash A\xra{}F(A,A)$ is an
isomorphism; this implies that we have $DA\Smash B=F(A,B)$ for all
$B$.  The category of dualisable objects is formally very similar to
the category of finite dimensional vector spaces over a field.

In Boardman's category $\CB$, or in the derived category $\CD_R$, it
is known that an object is dualisable iff it is small.  However,
$L_{K(n)}S^0$ is dualisable but not small in the $K(n)$-local
category.  A number of interesting things are known about
$K(n)$-locally dualisable spectra:
\begin{enumerate}
 \item A $K(n)$-local spectrum $X$ is dualisable iff
  $\dim_{K(n)_*}K(n)_*X<\infty$ (proved in~\cite{host:mkl}).
 \item If $X$ is a finite complex, then $L_{K(n)}\Sgi X$ is
  easily seen to be dualisable.
 \item If $X$ is a connected space with $|\prod_{k>0}\pi_kX|<\infty$,
  then it is probably true that $K(n)_*X$ is finite-dimensional and so
  $L_{K(n)}\Sgi X$ is dualisable.  The results in the literature
  involve some additional conditions, however; for example, the claim
  is true if $X=BG$ with $G$ finite~\cite{ra:mkf} (in which case
  $L_{K(n)}\Sgi BG$ turns out to be self-dual~\cite{st:kld}), or if
  $X$ is a double loop space~\cite{horawi:mha}.
\end{enumerate}

In~\cite{hopast:ash}, we assume that our generators are dualisable,
but not that they are small.  This theory has the advantage that any
localisation of a category satisfying our axioms, again satisfies our
axioms.  The disadvantages are:
\begin{itemize}
 \item[(i)] We need a smash product to formulate the definition of
  dualisability, and this is absent or unnatural in many examples,
  such as $\CD_R$ when $R$ is not commutative.
 \item[(ii)] If $\CC=\Ho(\CC_0)$ for some Quillen model category
  $\CC_0$, then there are natural conditions on objects in $\CC_0$
  guaranteeing that they are small in $\CC$.  This is not the case for
  dualisability: we need special geometric arguments to show that
  $G/H_+$ is dualisable in $\CB_G$, for example.
 \item[(iii)] There are naturally occurring cases where the generators
  are small but not dualisable, for example the category of
  $G$-spectra based on an incomplete universe, and possibly also
  derived categories for nonaffine schemes.
\end{itemize}

In~\cite{ne:tc}, Neeman introduces the notion of a
\emph{well-generated} triangulated category; the definition is
explained and simplified in~\cite{kr:nwg}.  To explain the nature of
this concept, we recall some generalisations of Quillen's small object
argument.  Quillen originally considered a category $\CE$ closed under
limits and colimits, and looked for objects $A$ that were small in the
sense that the functor $\CE(A,-)$ preserves filtered colimits.  Later,
it was realised that one can fix a large cardinal $\kp$ and say that
$A$ is \emph{$\kp$-small} if $\CE(A,-)$ preserves colimits of
sequences indexed by ordinals larger than $\kp$.  In many categories,
every object is $\kp$-small for some $\kp$, and in many applications
related to localisation, this is an adequate substitute for
smallness.  Neeman works with triangulated categories, which typically
do not have colimits for most diagrams.  Thus, the above cannot be
applied directly, but a somewhat more elaborate argument leads to
Neeman's well-generated categories.  It is shown in~\cite{ne:dcs} that
the derived category of any Grothendieck abelian category is
well-generated. 

\subsection{Representability}

A \emph{cohomology functor} on $\CC$ is a contravariant functor from
$\CC$ to the category $\Ab$ of abelian groups, that converts
coproducts to products and cofibre sequences to exact sequences.  For
any fixed object $Z$, it is well-known that the representable functor
$X\mapsto[X,Z]$ is a cohomology functor.  We say that \emph{the
  representability theorem holds for $\CC$} if the converse is true,
so that every cohomology theory on $\CC$ is representable.  Brown
proved the representability theorem for $\CB$, and the same proof
works for any compactly generated triangulated category.  By much more
elaborate arguments, Neeman has extended this to all well-generated
triangulated categories~\cite{ne:tc}, and similar results have been
obtained by Krause~\cite{kr:brt,kr:cfs} and Franke~\cite{fr:brt}.

It is also easy to see that if the representability theorem holds for
$\CC$, then it holds for any localisation of $\CC$.

In~\cite{hopast:ash}, we take the representability theorem as an
axiom; this gives a convenient way to treat compactly generated
categories and their localisations in parallel.  All other authors
assume axioms that turn out to imply the representability theorem.

\subsection{Extra axioms}

We now discuss some possible additional assumptions.  No author takes
any these as a standard axiom, but they define special classes of
examples with usefully simplified behaviour.  

\begin{itemize}
 \item[(a)] Let $\CC$ be an abelian category in which all
  monomorphisms split (and thus all epimorphisms split), and suppose
  we have an equivalence $\Sg\:\CC\xra{}\CC$.  We can give $\CC$ a
  triangulation by declaring that all all triangles of the form
  \[ A\op B \xra{f} B\op C \xra{g} C\op \Sg A \xra{h} \Sg A\op\Sg B
  \]
  (with $f(a,b)=(b,0)$ and so on) are cofibre sequences.  We call this
  the \emph{abelian case}.  Under mild finiteness conditions, one can
  show that $\CC$ is the category of graded $A_*$-modules, for some
  graded ring $A_*$ that is a finite product of graded division rings.
  This situation is discussed in~\cite[Section 8]{hopast:ash}.

  Examples include the category of rational $G$-spectra for any finite
  group $G$, or the category $\Stab_{kG}$ when $|G|$ is invertible in
  $k$. 

 \item[(b)] Suppose that $\CC$ has a symmetric monoidal structure, and
  that $\CC$ is generated by the single object $S$ (the unit for the
  smash product).  We then have a graded-commutative ring $\pi_*S$
  defined by $\pi_nS=[\Sg^nS,S]$.  If this is noetherian, we say that we
  are in the \emph{noetherian case}; this is discussed
  in~\cite[Section 6]{hopast:ash}.  Examples include $\CD_R$ where
  $R$ is commutative and noetherian, and $\CD_{kG}$.
  
 \item[(c)] Suppose again that $\CC$ has a symmetric monoidal
  structure, and that $\CC$ is generated by the single object $S$.  If
  $\pi_nS=0$ for $n<0$, then we are in the \emph{connective case};
  this is discussed in~\cite[Section 7]{hopast:ash}.
\end{itemize}

We will make a number of remarks about the noetherian case below.
Beyond that, we refer the reader to~\cite{hopast:ash} for further
discussion.

\section{Functors on small objects}
\label{sec-func}

Let $\CF$ be the category of small objects in $\CC$, and consider the
category $\CA=[\CF^{\text{op}},\Ab]$ of additive contravariant
functors from $\CF$ to the category of abelian groups.  This is a
bicomplete abelian category satisfying the AB5 condition (filtered
colimits are exact).  The functor $\Sg\:\CC\xra{}\CC$ induces a
functor $\Sg\:\CA\xra{}\CA$.  If $\CC$ has a good symmetric monoidal
structure, then so does $\CA$.  If the objects of $\CF$ are strongly
dualisable, then we have $\CF\simeq\CF^{\opp}$ and so
$\CA\simeq[\CF,\Ab]$.  One can think of $\CF$ as a ``ring with many
objects'', and regard $\CA$ as its module category.  

There is a Yoneda functor $h\:\CC\xra{}\CA$ sending $X$ to the functor
$h_X(A)=[A,X]$ (for $A\in\CF$).  The structure of $\CA$ and the
behaviour of $h$ have proved to be very useful in the study of $\CC$,
at least when $\CC$ is compactly generated.  If $\CC$ is merely
well-generated, then Neeman has developed a partially parallel theory
based on more complicated functor categories~\cite[Chapter 6]{ne:tc}.
In the compactly generated case, Beligiannis~\cite{be:rha} has
considered categories of the form $[\CG^{\opp},\Ab]$ where $\CG$ is an
arbitrary triangulated subcategory of $\CF$.

Let $\CE\sse\CA$ be the category of \emph{exact} functors: those that
send cofibre sequences in $\CF$ to exact sequences in $\Ab$.  It is
standard that $h_X\in\CE$ for all $X$.  In good cases, $\CE$ is the
category of objects of finite projective (or injective) dimension in
$\CA$, and $h\:\CC\xra{}\CE$ is close to being an equivalence; see
Section~\ref{sec-brown} for more discussion.

The functor $h\:\CC\xra{}\CA$ always preserves coproducts and sends
cofibre sequences to exact sequences.  In other words, it is an
$\CA$-valued homology theory on $\CC$.  A morphism $u\:X\xra{}Y$ in
$\CC$ is said to be \emph{phantom} if $h_u\:h_X\xra{}h_Y$ is zero.

\section{Types of subcategories}
\label{sec-types}

Let $R$ be a commutative noetherian ring, and let $\zar(R)$ denote the
space of prime ideals in $R$, with the Zariski topology.  (We do not
use the notation $\text{spec}(R)$, to avoid conflicting uses of the
word ``spectrum''.)  It turns out~\cite{ho:gmh,ne:ctd} that we can
recover $\zar(R)$ from the category $\CD_R$, in several slightly
different ways.  More precisely, one can recover the lattice of
radical ideals in $R$, which is well-known to be anti-isomorphic to
the lattice of closed subsets of $\zar(R)$, and this lattice
determines $\zar(R)$ itself~\cite{jo:ss}.  The key is to study various
lattices of subcategories of $\CD_R$.  Using parallel constructions in
the category $\Stab_{kG}$, we can recover the space
$\zar(H^*(G;\Fp))$.  In the case of Boardman's category, this study
makes contact with the chromatic approach to stable homotopy theory,
and the nilpotence theorems of Devinatz, Hopkins and Smith.  This has
many important applications that are not visible in the purely
algebraic examples.  For example, suppose we want to prove that all
finite spectra $X$ have some property $P(X)$.  Suppose we can show
that
\begin{itemize}
 \item Whenever we have a cofibre sequence $X\xra{}Y\xra{}Z$ in which two
  terms have property $P$, then the third also has property $P$
 \item Whenever $P(X\Wedge Y)$ holds, so do $P(X)$ and $P(Y)$
 \item There exists a finite spectrum $X$ such that $H_*(X;\Q)\neq 0$
  and $P(X)$ holds.
\end{itemize}
Then one can show using the subcategory classification theorems that
$P(X)$ holds for all $X$.

The basic definitions are as follows.
\begin{definition}
 Let $\CD$ be a full subcategory of $\CC$.  For simplicity, we assume
 that any object in $\CC$ that is isomorphic to an object in $\CD$, is
 itself in $\CD$.  Let $\CA$ be an arbitrary collection (possibly a
 proper class) of objects in $\CC$.
 \begin{itemize}
  \item[(a)] $\CD$ is \emph{thick} if 
   \begin{itemize}
    \item[(i)] The zero object lies in $\CD$.
    \item[(ii)] Any retract of any object in $\CD$, again lies in
     $\CD$.
    \item[(iii)] Whenever $X\xra{}Y\xra{}Z$ is a cofibre sequence with two
     terms in $\CD$, the third term is also in $\CD$.
   \end{itemize}
  \item[(b)] If $\CC$ has a symmetric monoidal structure, we say that
   $\CD$ is an \emph{ideal} if $X\Smash Y\in\CD$ whenever $X\in\CD$.
   Dually, we say that $\CD$ is a \emph{coideal} if $F(Y,Z)\in\CD$
   whenever $Z\in\CD$.
  \item[(c)] $\CD$ is a \emph{localising subcategory} if it is thick,
   and closed under (possibly infinite) coproducts.  Dually,
   $\CD$ is a \emph{colocalising subcategory} if it is thick,
   and closed under (possibly infinite) products.
  \item[(d)] A \emph{(co)localising (co)ideal} is a (co)localising
   subcategory that is also a (co)ideal.
  \item[(e)] A \emph{bilocalising subcategory} is a subcategory that
   is both a localising subcategory and a colocalising subcategory.  A
   \emph{biideal} is a subcategory that is both an ideal and a
   coideal.
 \end{itemize}
\end{definition}

If $\CC$ is monoidal and the unit object $S\in\CC$ is small and
generates $\CC$, then every (co)localising subcategory is a (co)ideal.
This holds in the following cases:
\begin{itemize}
 \item[(1)] $\CC=\CB$ (but not $\CC=\CB_G$ for general $G$)
 \item[(2)] $\CC=\CD_R$, where $R$ is commutative
 \item[(3)] $\CC=\Stab_{kG}$, where $k$ has characteristic $p$ and $G$
  is a $p$-group.
\end{itemize}

If we have a functor $F$ between triangulated categories that
preserves cofibre sequences, then $\ker(F):=\{X\st FX=0\}$ is evidently a
thick subcategory.  Similarly, if $F$ is a functor from a triangulated
category to an abelian category, and $F$ converts cofibre sequences to
exact sequences, then $\ker(F)$ will again be a thick subcategory.
Under various auxiliary conditions, we can conclude that $\ker(F)$ is
a (co)localising subcategory or a (co)localising ideal.

Little is known about classification of subcategories that are not
ideals.  The examples studied in~\cite[Section 6]{becari:tss} suggest
that there is no simple and general picture.

On the other hand, there are good classification results for many of
our central examples; the main method of proof will be discussed in
Section~\ref{sec-nilpotence}.  To state the results, it is convenient
to introduce one more definition.  Given an object $A\in\CF$, we write
$\thickid{A}$ for the smallest thick ideal containing $A$.  We then
say that a thick ideal $\CI$ is \emph{finitely generated} if it has
the form $\thickid{A}$ for some $A$ (this makes sense because the
thick ideal generated by $A_1,\ldots,A_r$ is also generated by the
single object $A=\bigWedge_iA_i$).  A classification of finitely
generated thick ideals (or thick subcategories) in $\CF$ extends in a
fairly obvious way to give a classification of all ideals (or
thick subcategories) in $\CF$.
\begin{itemize}
 \item[(a)] Let $R$ be a commutative noetherian ring, and put
  $\CC=\CD_R$.  Then the localising subcategories of $\CC$ biject with
  the colocalising subcategories, and with the subsets of $\zar(R)$.
  The finitely generated thick subcategories of $\CF$ biject with
  closed subsets of $\zar(R)$.  On the other hand, Neeman has
  considered the nonnoetherian ring 
  \[ R = k[x_2,x_3,x_4,\ldots]/(x_2^2,x_3^3,x_4^4,\ldots), \]
  where $k$ is a field.  This has only one prime ideal, but $\CD_R$
  has an enormous collection of localising
  subcategories~\cite{ne:obc}.
 \item[(b)] Let $k$ be a field, let $G$ be a finite group, and put
  $\CC=\Stab_{kG}$.  It is proved in~\cite{becari:tss} that the
  finitely generated thick ideals in $\CF$ biject with the closed
  subsets of the projective scheme $\proj(H^*(G;k))$.  One can also
  show that the (co)localising (co)ideals in $\CC$ biject with all
  subsets of $\proj(H^*(G;k))$.
 \item[(c)] In the category of $E(n)$-local spectra, the
  (co)localising subcategories in $\CC$ biject with subsets of
  $\{0,1,\ldots,n\}$, and the thick subcategories of $\CF$ biject with
  the subsets of the form $\{m,m+1,\ldots,n\}$ for some $m$.  All the
  relevant subcategories are (co)ideals~\cite[Theorem 6.14]{host:mkl}.
 \item[(d)] Now let $\CC$ be the category of $K(n)$-local spectra.
  Then $0$ and $\CC$ are the only localising subcategories of $\CC$,
  and also the only colocalising subcategories of $\CC$~\cite[Theorem
  7.5]{host:mkl}.  Similarly, $0$ and $\CF$ are the only thick
  subcategories of $\CF$.
 \item[(e)] Finally, let $\CC$ be the category of $p$-local spectra.
  The thick subcategories of $\CF$ are then the categories
  $\CF_n:=\{X\st K(m)_*X=0 \text{ for all } m<n\}$, where
  $0\leq n\leq\infty$; this was proved in~\cite{hosm:nshii}.  The
  theory of Bousfield classes gives many known examples of
  (co)localising subcategories and inclusions between them.  However,
  almost nothing is known about the collection of \emph{all}
  localising subcategories (which might even be a proper class).
\end{itemize}
The classification results in Examples~(a) and~(b) have a partial
generalisation that applies in the noetherian case.  A strong but
technically complex statement is proved in~\cite[Section
6.3]{hopast:ash}; given some additional hypotheses (conjecturally
always satisfied) this implies the evident analog of~(a) and~(b).

So far we have only discussed results about ideals in $\CF$; we next
consider results about (co)localising subcategories or (co)ideals in
$\CC$.  On the one hand, given $\CD\sse\CC$ we can certainly consider
the thick subcategory $\CD\cap\CF\sse\CF$; if we have a good
understanding of $\CF$ then this will be a useful invariant, but
rather a coarse one.  On the other hand, given a thick subcategory
$\CA\sse\CF$ we can consider the category 
\[ \CA^\perp := \{X\st [A,X]_*=0 \text{ for all } A\in\CA\}, \]
which is easily seen to be a bilocalising subcategory.  The
\emph{telescope conjecture} for $\CC$ is closely related to the
statement that every bilocalising subcategory is of the form
$\CA^\perp$ for some $\CA$.  This is known to hold in $\CD_R$ when $R$
is noetherian and commutative, and also in $\Stab_{kG}$ when $k$ has
characteristic $p$ and $G$ is a $p$-group.  It is believed to be false
in Boardman's category, although many years of study have still not
produced a watertight argument.

\section{Quotient categories and Bousfield localisation}
\label{sec-quotient}

Let $\CC$ be a triangulated category.  Given a thick subcategory
$\CD$, we can look for a triangulated category $\CC'$ and an exact
functor $Q\:\CC\xra{}\CC'$ that sends all objects in $\CD$ to zero.
It turns out that there is an initial example of such a functor, whose
target we call $\CC/\CD$.  To be more precise, we say that a map
$s\:X\xra{}Y$ in $\CC$ is a $\CD$-equivalence if the cofibre of $s$
lies in $\CD$.  The class of $\CD$-equivalences has a number of useful
properties:
\begin{itemize}
 \item Any isomorphism is a $\CD$-equivalence.
 \item Given morphisms $X\xra{s}Y\xra{t}Y$, if any two of $\{s,t,ts\}$
  are $\CD$-equivalences then so is the third.
 \item Given maps $X\xra{f}Y\xla{s}Z$ in which $s$ is a
  $\CD$-equivalence, there is a commutative square
  \[ \xymatrix{
      W \dto_t \rto^g & Z \dto^s \\
      X \rto_f & Y
     }
  \]
  in which $t$ is a $\CD$-equivalence.
\end{itemize}
We then define $\CC/\CD$ as follows: the objects are the same as in
$\CC$, and the morphisms from $X$ to $Y$ are equivalence classes of
``formal fractions'' $gt^{-1}$, where $g$ and $t$ fit in a diagram of
the shape $X\xla{t}W\xra{g}Y$, and $t$ is a $\CD$-equivalence.  The
properties listed above allow us to compose and manipulate fractions
in a natural way.

Krause has considered some more delicate notions of quotient
categories, which are important in the study of smashing
localisations; but we will not discuss these here.

As is well-known, there is a potential problem with the above
construction, which is of great importance in some applications.  We
always assume implicitly that the morphism sets $\CC(X,Y)$ are genuine
sets; but as defined above, $(\CC/\CD)(X,Y)$ might be a proper class.
There are a number of techniques that can be used in different
circumstances to show that this problem does not arise.  As far as I
know, no one has looked systematically for examples where proper
classes do arise; it is possible that (some version of) our standing
axioms are enough to prevent this.

If $\CC/\CD$ has small Hom sets, then for any $Y\in\CC$ we have a
functor from $\CC$ to $\Ab$ given by $X\mapsto(\CC/\CD)(X,Y)$.  The
Representability Theorem shows that this is representable, so we have
an object $LX\in\CC$ and an isomorphism
$\CC(X,LY)\simeq(\CC/\CD)(X,Y)$, naturally in $X$.  A standard
argument shows that $L$ can be regarded as a functor
$\CC/\CD\xra{}\CC$, right adjoint to the quotient functor
$\CC\xra{}\CC/\CD$.  If we put 
\[ \CD^\perp = \{Y \st \CC(X,Y)=0 \text{ for all } X\in\CD \}, \]
we find that $L$ actually gives an equivalence
$\CC/\CD\simeq\CD^\perp$.  We also use the letter $L$ for the
composite functor $\CC\xra{}\CC/\CD\xra{}\CD^\perp\sse\CC$; in this
guise, it is left adjoint to the inclusion of $\CD^\perp$ in $\CC$.  

The functors $L\:\CC\xra{}\CC$ arising in this way can be
characterised by certain well-known properties: they are exact
functors, equipped with a natural map $i_X\:X\xra{}LX$ such that
$Li_X\:LX\xra{}L^2X$ is an equivalence, and
$i_X^*\:[LX,LY]\xra{}[X,LY]$ is an isomorphism for all $Y$.  We call
such a pair $(L,i)$ a \emph{Bousfield localisation functor}.  We can
recover $\CD$ as the category $\ker(L)=\{X\st LX=0\}$.  

The above discussion shows that quotients are really localisations.
Of course, the converse is also true: to invert a class of maps $\CE$
is the same as to quotient out the localising subcategory generated by
the cofibres of the maps in $\CE$.

In~\cite{hopast:ash}, it is assumed that $\CC$ is symmetric monoidal
and that $\CD$ is a localising ideal.  In this case, the quotient
category $\CC/\CD$ (or equivalently, the category $\CD^\perp$)
inherits a symmetric monoidal structure.

Given any localisation functor $L$, there is another functor $C$ and
natural transformations $CX\xra{q_X}X\xra{i_X}LX\xra{d_X}\Sg CX$
giving a cofibre sequence for all $X$.  The theory can be set up in
such a way that $C$ and $L$ play precisely dual r\^oles.  In any case,
the pair $(L,i)$ determines $(C,q)$ (up to an obvious notion of
equivalence) and \emph{vice versa}.

Given two localisation functors $(L,i)$ and $(L',i')$, there is at
most one morphism $u\:L\xra{}L'$ with $ui=i'$.  We write $L\geq L'$ if
such a morphism exists.  This gives a partial order on the collection
of isomorphism classes of localisation functors.  It is not known
whether this collection is a set or a proper class.  

\section{Versions of the Bousfield lattice}
\label{sec-lattice}

In this section, we assume that $\CC$ has a symmetric monoidal
structure.  Without such a structure, one could set up a formal theory
along the same lines, but it seems hard to analyze any examples
explicitly.

We can now define various partially ordered sets $\Lm_i$; some of them
may actually be proper classes, but we will suppress this from the
terminology.  An optimistic conjecture would be that they are all the
same; this is known to be true in the noetherian case.  The general
case appears to be open (related work of Gutierrez and Casacuberta
turns out not to provide a counterexample).

\begin{definition}
 \begin{itemize}
  \item $\Lm_0$ is the class of all colocalising coideals (ordered by
   inclusion).  For any class $\CA$ of objects in $\CC$, we write
   $\colocid{\CA}$ for the intersection of all colocalising ideals
   containing $\CA$.  The poset $\Lm_0$ is actually a lattice, with
   meet operation $\CD\wedge\CD'=\CD\cap\CD'$, and join
   $\CD\vee\CD'=\colocid{\CD\cup\CD'}$.
  \item $\Lm_1$ is the class of all localising ideals, ordered by
   \emph{reverse} inclusion.  This is a lattice by a dual argument.
  \item For any class $\CA$ of objects in $\CC$, we put 
   \begin{align*}
    \CA^\perp   &= \{X \st F(A,X) = 0 \text{ for all } A\in\CA\} \\
    {}^\perp\CA &= \{X \st F(X,A) = 0 \text{ for all } A\in\CA\}.
   \end{align*}
   It is easy to see that $\CA^\perp\in\Lm_0$ and
   ${}^\perp\CA\in\Lm_1$.  We say that a colocalising coideal $\CD$ is
   \emph{closed} if it has the form $\CA^\perp$ for some $\CA$, or
   equivalently if $\CD=({}^\perp\CD)^\perp$; we write $\Lm_2$ for
   the set of closed colocalising coideals, so $\Lm_2\sse\Lm_0$.  
  \item Dually, we say that a localising ideal $\CE$ is \emph{closed}
   if it has the form ${}^\perp\CA$ for some $\CA$, or equivalently if
   $\CE={}^\perp(\CE^\perp)$.  We write $\Lm_3$ for the set of closed
   localising ideals, so that $\Lm_3\sse\Lm_1$.
  \item There are order-preserving maps $\Lm_0\xra{}\Lm_3$ and
   $\Lm_1\xra{}\Lm_2$, given by $\CD\mapsto{}^\perp\CD$ and
   $\CE\mapsto\CE^\perp$.  A purely formal argument shows that these
   give an isomorphism $\Lm_2\simeq\Lm_3$.
  \item We say that a colocalising coideal $\CD$ is \emph{reflective}
   if the inclusion $\CD\xra{}\CC$ has a left adjoint; one can show
   that this implies that $\CD$ is closed, so the coreflective
   coideals give a subset $\Lm_4\sse\Lm_2$.  Dually, we say
   that a localising ideal $\CE$ is \emph{coreflective} if the
   inclusion has a right adjoint, and these ideals give a subset
   $\Lm_5\sse\Lm_3$.  The bijection $\Lm_2\simeq\Lm_3$ restricts to
   give a bijection $\Lm_4\simeq\Lm_5$.  If $\CD$ and $\CE$ correspond
   under this bijection, then there is a pair of functors $(L,C)$ as
   in Section~\ref{sec-quotient}, with
   \begin{align*}
    \CD &= \img(L) = \ker(C) = \ker(L)^\perp = \img(C)^\perp \\
    \CE &= \img(C) = \ker(L) = {}^\perp\ker(L) = {}^\perp\img(C).
   \end{align*}
   It follows that $\Lm_4$ and $\Lm_5$ are equivalent to the poset of
   localisation functors $L$ for which $\ker(L)$ is an ideal, or to
   the poset of colocalisation functors $C$ for which $\ker(C)$ is a
   coideal. 
  \item We say that a localising ideal $\CE$ is \emph{principal} if
   $\CE=\locid{\{E\}}$ for some object $E$.  Note that if
   $\CE=\locid{\{E_i\}_{i\in I}}$ (where $I$ is a set, not a proper
   class) then we also have $\CE=\locid{\{E\}}$ where
   $E=\bigWedge_{i\in I}E_i$, so $\CE$ is principal.  In Boardman's
   category, it is known that principal ideals are coreflective, so
   they form a subset $\Lm_6\sse\Lm_5$.  It is not clear in what
   generality this argument works.  If $\CE=\locid{\{E\}}$ then the
   corresponding localisation functor is called \emph{stable
     $E$-nullification}, and written $P_{\Sg^*E}$.  (Confusingly, it
   was called \emph{colocalisation} in Bousfield's original papers.)
  \item We say that a localising ideal $\CE$ is a \emph{Bousfield
     class} if it has the form
   $\CE=\bc{E}=\{X\st E\Smash X=0\}$ for some $E\in\CC$.  We
   write $\Lm_7$ for the collection of all Bousfield classes.  In
   Boardman's category, this is contained in $\Lm_6$; it is not clear
   how far this fact can be generalised.  It is also known that
   $\Lm_7$ is a set rather than a proper class~\cite{oh:ihh,dwpa:ot}.
   To see this, for any finite spectrum $A$ and any element
   $x\in E_*A$ we let $\ann_A(x)$ be the set of maps $f\:A\xra{}B$ in
   $\CF$ such that $(E_*f)(x)=0$.  We then write
   \[ \langle\langle E\rangle\rangle = 
        \{\ann_A(x)\st A\in\CF\;,\; x\in E_*A\},
   \]
   and call this the \emph{Ohkawa class} of $E$.  As $\CF$ has small
   Hom sets and only a set of isomorphism classes, we see that there
   is only a set of possible Ohkawa classes.  One can check that
   $\langle\langle E\rangle\rangle$ determines $\bc{E}$, so
   there is only a set of Bousfield classes.
 \end{itemize}
\end{definition}

\section{Special types of localisation}
\label{sec-special}

In this section, we assume that $\CC$ is compactly generated.

\begin{definition}
 A Bousfield localisation functor $L\:\CC\xra{}\CC$ is \emph{smashing}
 if $\img(L)$ (which is automatically a colocalising subcategory) is
 closed under coproducts (and so is also a localising subcategory).
 In the monoidal case, this implies that there is a natural
 equivalence $LS\Smash X\xra{}LX$, and the corresponding
 colocalisation functor $C$ also satisfies $CX=CS\Smash X$.  The
 category $\CD=\img(L)=\ker(C)=\ker(L)^\perp=\img(C)^\perp$ is then
 both a localising ideal and a colocalising coideal.  We put
 $\CU={}^\perp\CD$ and $\CV=\CD^\perp$, and then $\hC X=F(LS,X)$ and
 $\hL X=F(CS,X)$.  It turns out that $\hL$ is a localisation functor,
 and $\hC$ is the corresponding colocalisation.  Moreover, we have
 \begin{align*}
  \CU &= \img(C) = \ker(L) = \bc{LS} \\
  \CV &= \img(\hL) = \ker(\hC) \\
  \CD &= \img(\hC) = \ker(\hL) = \img(L) = \ker(C) \\
      &= \CU^\perp = {}^\perp\CV = \bc{CS}.
 \end{align*}
 It follows that $C\hC=0=\hL L$, so $\hL C\simeq\hL$ and
 $C\hL\simeq C$.  This implies that $\hL\:\CU\xra{}\CV$ and
 $C\:\CV\xra{}\CU$ are mutually inverse equivalences.  
 
 Apart from the finite localisations discussed below, the most
 important examples are the localisations with respect to the
 Johnson-Wilson spectra $E(n)$.  It is a highly nontrivial
 theorem~\cite[Chapter 8]{ra:nps} that these are smashing.

 Krause~\cite{kr:sst} has shown that $L$ is determined by the set
 $\ann(L|_{\CF})$ of morphisms $u\:A\xra{}B$ in $\CF$ for which
 $Lu=0$.  This means in particular that there is only a set of
 smashing localisations.
\end{definition}

\begin{definition}
 A \emph{finite localisation} is a localisation functor
 $L\:\CC\xra{}\CC$ where $\ker(L)=\loc{\CA}$ for some thick
 subcategory $\CA\sse\CF$.  Functors of this type are always
 smashing~\cite{mi:fl}\cite[Section 3.3]{hopast:ash}.  One formulation
 of the \emph{telescope conjecture} for $\CC$ is the statement that
 every smashing localisation is a finite localisation.  This is known
 to be true in many noetherian cases, but believed to be false in
 Boardman's category.  Keller~\cite{ke:gsc} has provided a
 counterexample in $\CD_R$ for a certain ring $R$ (but his framework
 of definitions is slightly different from ours, and we have not
 pinned down the precise relationship).
 
 An important example of finite localisation is as follows.  Let $R$
 be a noetherian ring, and put $\CC=\CD_R$.  Fix an ideal $I$, and let
 $\CA$ consist of the objects $X\in\CC$ for which $\pi_*X$ is an
 $I$-torsion module.  Here the category $\CU=\ker(L)$ consists of the
 $I$-torsion objects in $\CC$, and $\CV=\ker(\hC)$ consists of
 $I$-complete objects in a suitable sense.  Thus, the equivalence
 $\CU\simeq\CV$ shows that the torsion category and the complete
 category are essentially the same.  All this is closely related to
 the theory of local (co)homology~\cite{grma:dfi}.
 See~\cite{dwgr:cmt} and~\cite[Section 3.3]{hopast:ash} for other
 perspectives.
\end{definition}

\begin{definition}
 Suppose that $\CC$ is symmetric monoidal, and is generated by the
 unit object $S$.  Given a set $T$ of homogeneous elements in the
 graded ring $\pi_*S$, we let $\CA$ denote the thick subcategory of
 $\CF$ generated by the cofibres of the maps in $T$, and then let $L$
 be the corresponding finite localisation functor.  One can show that
 $\pi_*LX=(\pi_*X)[T^{-1}]$.  Functors of this type are called
 \emph{algebraic localisations}; in the special case where $T\sse\Z$,
 they are called \emph{arithmetic localisations}.  
\end{definition}

\section{Nilpotence}
\label{sec-nilpotence}

Our understanding of Boardman's category relies heavily on the
nilpotence theorem of Devinatz, Hopkins and Smith~\cite{dehosm:nsh}
and its consequences~\cite{hosm:nshii,ra:nps}.  We next explain the
formal parts of this story that are amenable to axiomatic
generalisation~\cite{st:ebc}\cite[Section 5]{hopast:ash}.  We will
assume here that $\CC$ is compactly generated and has a symmetric
monoidal structure, and that all objects of $\CF$ are strongly
dualisable.

We say that an object $I\in\CF$ equipped with a map $i\:I\xra{}S$ is
an \emph{ideal} if the map $i\Smash 1\:I\Smash S/I\xra{}S/I$ is null.
(Here $S/I$ denotes the cofibre of $i$.)  We write $I\leq J$ if the
map $I\xra{}S\xra{}S/J$ is zero.  It turns out that if $I\xra{i}S$ and
$J\xra{j}S$ are ideals, then so is $I\Smash J\xra{i\Smash j}S$; we
will just write $IJ$ for this, making the set of isomorphism classes
of ideals into a commutative monoid.  We say that $I$ and $J$ are
\emph{radically equivalent} if for large $n$ we have $I^n\leq J$ and
$J^n\leq I$.  We write $\bId(S)$ for the set of radical equivalence
classes of ideals.  Given any $A\in\CF$, the fibre of the unit map
$S\xra{}F(A,A)$ is an ideal, and we write $\ann(A)$ for its
equivalence class.  One can show that the rule
$\ann(A)\leftrightarrow\thickid{A}$ gives a well-defined bijection
between $\bId(S)$ and the set of finitely generated thick ideals in
$\CF$.

Next, we say that a map $u\:A\xra{}B$ in $\CF$ is
\emph{smash-nilpotent} if the $m$'th smash power
$u^{(m)}\:A^{(m)}\xra{}B^{(m)}$ is zero for $m\gg 0$.  One checks that
$I^m\leq J$ for some $m$ iff the composite $I\xra{}S\xra{}S/J$ is
smash-nilpotent.  

Now suppose we are given a set $N$ and a collection of objects
$K(n)\in\CC$ for each $n\in N$.  For any object $X\in\CC$ we put
$\supp(X)=\{n\st K(n)\Smash X\neq 0\}$.  Similarly, given a thick
ideal $\CA\sse\CF$ we put $\supp(\CA)=\bigcup_{A\in\CA}\supp(A)$.

We say that the $K(n)$'s \emph{detect ideals} if whenever $A,B\in\CF$
and $\supp(A)\sse\supp(B)$, we have $\thickid{A}\sse\thickid{B}$.
This implies that the map $\CA\xra{}\supp(\CA)$ gives an embedding of
the lattice of thick ideals in the lattice of subsets of $N$.  (Except
in the noetherian case, we know of no general method to determine the
image of this map.)

We next explain two versions of what it might mean for the $K(n)$'s to
\emph{detect nilpotence}.  It is a key theorem that if the $K(n)$'s
detect nilpotence, then they also detect ideals.

For the most algebraically natural version, we need auxiliary
hypotheses.  First, we assume that each $K(n)$ has a commutative ring
structure, and that every nonzero homogeneous element in the
coefficient ring $K(n)_*$ is invertible (so $K(n)_*$ is a graded
field).  We also assume that the resulting K\"unneth maps
\[ K(n)_*(X)\ot_{K(n)_*}K(n)_*(Y) \xra{} K(n)_*(X\Smash Y) \]
are isomorphisms for all $X$ and $Y$ (this is not automatic unless
$\CC$ is generated by $\{S\}$).  Thus, we can regard $K(n)$ as giving
a monoidal functor from $\CF$ to the category of finite-dimensional
vector spaces over $K(n)_*$.  We say that a map $u\:A\xra{}B$ in $\CF$
is $K(*)$-null if the induced map $K(n)_*A\xra{}K(n)_*B$ is zero for
all $n$.  We say that the $K(n)$'s \emph{detect smash-nilpotence} if
every $K(*)$-null map is smash-nilpotent.  Assuming this, a
straightforward argument (based on our discussion of $\bId(S)$) shows
that the $K(n)$'s detect ideals.  

In a general stable homotopy category $\CC$, it is very hard to
produce ring objects $K(n)$ such that $K(n)_*$ is a graded field.
However, one can use another line of argument with rather different
hypotheses.  First, we say that the $K(n)$'s \emph{detect rings} if
for every nonzero ring object $R$ we have $K(n)\Smash R\neq 0$ for
some $n$.  (This will obviously hold if 
$\bigWedge_n\bc{K(n)}=\bc{S}$ in the Bousfield
lattice.)  Suppose in addition that whenever $A\in\CF$ and 
$K(n)\Smash A\neq 0$ we have 
$\bc{K(n)\Smash A}=\bc{K(n)}$.  We claim that the
$K(n)$'s detect ideals.  To see this, consider a thick ideal
$\CA\sse\CF$, and let $L$ be the finite localisation functor with
$\ker(L)=\loc{\CA}$ (and thus $\ker(L)\cap\CF=\CA$).  Given $X\in\CF$
with $\supp(X)\sse\supp(\CA)$, we must show that $X\in\CA$.  It turns
out to be equivalent to say that the ring object $R=F(X,X)\Smash LS$
is zero, so it will suffice to show that $K(n)\Smash R=0$ for all $n$,
and this is easy.

The main examples are as follows.
\begin{itemize}
 \item[(a)] In the motivating
  example~\cite{dehosm:nsh,hosm:nshii,ra:nps}, $\CC$ is the category of
  $p$-local spectra, $N$ is $\N\cup\{\infty\}$, and $K(n)$ is the
  $n$'th Morava $K$-theory (which is well-known to be a field theory,
  up to a slight adjustment of the definition at the prime $2$).  The
  proof that these theories detect nilpotence is a \emph{tour de force}
  of stable homotopy theory, using methods very far from those surveyed
  in this paper.  It follows that they also detect smash nilpotence.
  It is deduced in~\cite{hosm:nshii} that they detect nilpotence in
  various other senses, and that they detect rings.  Part of this
  argument can be axiomatised (at least in a connective stable homotopy
  category) but we shall not attempt that here.
 \item[(b)] Let $G$ be a finite group, and let $\CC$ be the category
  of $p$-local $G$-spectra.  We then let $N$ be the set of pairs
  $(H,n)$, where $H$ is a (representative of a) conjugacy class of
  subgroups of $G$, and $n\in\N\cup\{\infty\}$.  We take $K(H,n)$ to
  be the representing object for the cohomology theory
  $X\mapsto K(n)^*\Phi^HX$, where $\Phi^H$ is the geometric fixed
  point functor, and $K(n)$ is the usual nonequivariant Morava
  $K$-theory.  These representing objects can be made quite explicit,
  but we shall not give the details.  It follows quite easily from the
  previous example that they detect smash-nilpotence, and thus that
  they detect ideals~\cite{st:ebc}.
 \item[(c)] Let $R$ be a noetherian ring, and put $\CC=\CD_R$.  We
  then take $N$ to be the set of prime ideals in $R$, and let
  $K(\pri)$ be the field of fractions of $R/\pri$.
  Here it is not hard to show that
  $\bigWedge_{\pri}\bc{K(\pri)}=\bc{S}$ and that $\bc{K(\pri)}$ is a
  minimal Bousfield class; it follows that these objects detect
  nilpotence, and also that they detect
  ideals~\cite{ne:ctd}\cite[Section 6]{hopast:ash}. 
 \item[(d)] Now consider the case $\CC=\Stab_{kG}$, where $k$ is a
  field of characteristic $p$ and $G$ is a finite $p$-group.  Take $N$
  to be the set of homogeneous prime ideals in $H^*(G;k)$.  Next, fix
  an algebraically closed field $L$ of infinite transcendence degree
  over $k$.  For any $\pri\in N$, the theory of ``shifted subgroups''
  gives an algebra $A\leq LG$ isomorphic to $L[u]/u^p$ and an object
  $K(\pri)\in\CC$ such that $K(\pri)\Smash M=0$ iff $L\ot_kM$ is free
  as a module over $A$.  It follows easily from the infinite version
  of Dade's Lemma~\cite{becari:cviii} that
  $\bigWedge_{\pri}\bc{K(\pri)}=\bc{S}$.  Morever, we see
  from~\cite[Theorem 10.8]{becari:cviii} that
  $\bc{K(\pri)\Smash M}=\bc{K(\pri)}$ whenever $K(\pri)\Smash M\neq
  0$, so the $K(\pri)$'s detect ideals.
\end{itemize}

\section{Brown representability}
\label{sec-brown}

For all authors, it is either an axiom or a theorem that cohomology
functors defined on the whole category $\CC$ are representable.  It
follows easily that the Yoneda functor is an equivalence between $\CC$
and the category of cohomology functors defined on $\CC$.  This is a
very satisfactory result, with many applications (existence of
infinite products, existence of Bousfield localisations,
Brown-Comenetz duality, and so on).

It is desirable to extend this result to various subcategories
$\CD\sse\CC$.  If there is an exact localisation functor
$L\:\CC\xra{}\CD$ (as in Section~\ref{sec-quotient}), then this is
easy.  In some other cases, it can be proved using Neeman's theory of
well-generated categories~\cite[Chapter 8]{ne:tc}.

Similarly, it would be helpful to have a dual theorem.  This should
say that any product-preserving exact covariant functor $\CC\xra{}\Ab$
has the form $Y\mapsto[X,Y]$ for some representing object $X$.  This
has also been proved by Neeman~\cite{ne:brd,ne:tc}, under some
additional hypotheses.

Next, let $\CF\sse\CC$ be the subcategory of small objects, and
suppose that $\CF$ generates $\CC$.  Given a cohomology functor
$H\:\CF^{\opp}\xra{}\Ab$, it is natural to ask whether there is an
object $Z\in\CC$ and a natural isomorphism $HX=[X,Z]$ for $X\in\CF$.
It is equivalent to ask whether $H$ can be extended to a cohomology
functor defined on all of $\CC$.

We first observe that in the case $\CC=\CB$, this reduces to a more
familiar question.  In that context, the Spanier-Whitehead duality
functor $D\:X\mapsto F(X,S)$ gives an equivalence
$\CF^{\opp}\simeq\CF$, so the covariant functor 
$H'=H\circ D\:\CF\xra{}\Ab$ is homological.  A natural isomorphism
$HX=[X,Z]$ (for all $X\in\CF$) is thus the same as a 
$H'X=\pi_0(X\Smash Z)$, and Brown's homological representability
theorem (in the version proved by Adams) says that such an isomorphism
can always be found.  Adams's proof used some countability arguments,
and implicitly relied on the existence of an underlying model
category, so it could not directly be transferred to our axiomatic
setting.  Margolis~\cite{ma:ssa} and Neeman~\cite{ne:tba}
independently gave reformulations that do not use model categories.
Neeman also showed, however, that the countability hypothesis is
essential. 

To explain this and related results (mostly distilled
from~\cite{ne:tba,chkene:fbr,be:rha}), it is convenient to use the
category $\CA=[\CF^{\opp},\Ab]$, the subcategory $\CE\sse\CA$ of exact
functors, and the Yoneda functor $h\:\CC\xra{}\CE$, as discussed in
Section~\ref{sec-func}.  Brown's theorem says that when $\CC=\CB$, the
functor $h\:\CC\xra{}\CE$ is full and essentially surjective.  The
work of Margolis and Neeman says that the same holds whenever
\begin{itemize}
 \item[(a)] $\CC$ is compactly generated; and
 \item[(b)] $\CF$ has only countably many isomorphism classes, and
  $\CF(A,B)$ is countable for all $A,B\in\CF$. 
\end{itemize}

Now consider the case $\CC=D(k[x,y])$, where $k$ is a field.
Neeman has shown that if $|k|\geq\aleph_2$ then $h$ is not full, and
if $|k|\geq\aleph_3$ then $h$ is not essentially surjective.  These
examples are obtained from more general and more complicated
statements of two different types.  Firstly, there are results
relating properties of $h$ to homological algebra in $\CA$; secondly,
there are relations between homological algebra in $\CA$ and in the
category $\CM_R$ of $R$-modules, in the case where $\CC=\CD_R$.

For the first step, we define $\pgldim(\CC)$ to be the supremum of the
projective dimensions in $\CA$ of all the objects in $\CE$.  Even
though $h$ need not be essentially surjective, this is known to be the
same as the supremum of the projective dimensions of objects in the
image of $h$.  It is also known that
\begin{itemize}
 \item[(a)] $\pgldim(\CC)\leq 1$ iff $h$ is full.
 \item[(b)] If $\pgldim(\CC)\leq 2$, then $h$ is essentially
  surjective.
 \item[(c)] Thus, if $h$ is full, then it is essentially surjective.
\end{itemize}

For the second step, we recall some additional definitions.  An
$R$-module $P$ is said to be \emph{pure-projective} if it is a retract
of a (possibly infinite) direct sum of finitely presented modules.
The \emph{pure-projective dimension} of a module $M$ is the minimum
possible length for a pure projective resolution of $M$.  The
\emph{pure global dimension} of $R$ (written $\pgldim(R)$) is the
supremum of the pure-projective dimensions of all $R$-modules.  The
ring $R$ is said to be \emph{hereditary} if every submodule of a
projective module is again projective.  It is known that
$\pgldim(R)\leq\pgldim(\CD_R)$; the inequality is an equality when $R$
is hereditary, but can be strict in more general cases.  Using this
and some additional arguments, one proves the following result:
\begin{theorem}
 Suppose that $\CC=\CD_R$, where $R$ is hereditary.  Then the functor
 $h\:\CC\xra{}\CA$ is
 \begin{itemize}
  \item[(a)] full iff $\pgldim(R)\leq 1$
  \item[(b)] essentially surjective iff $\pgldim(R)\leq 2$.
 \end{itemize}
\end{theorem}

Benson and Gnacadja~\cite{begn:pmpi,begn:pmpii} have proved similar
results for the case $\CC=\Stab_{kG}$, involving questions of purity
for $kG$-modules.  In particular, they show that the following are
equivalent:
\begin{itemize}
 \item[(a)] $h$ is full and essentially surjective
 \item[(b)] $\pgldim(kG)\leq 1$
 \item[(c)] Either $k$ is countable, or the Sylow $p$-subgroup of $G$
  is cyclic (where $p$ is the characteristic of $k$). 
\end{itemize}
They also give a number of intricate examples related to these
results.

Now consider the case where $h$ is full and essentially surjective, as
with the original case of Boardman's category of spectra.  We then say
that $\CC$ is a \emph{Brown category}.  This has a number of useful
consequences~\cite{hopast:ash,chst:pmh,be:rha}.  Firstly, for
$F\in\CA$, the following are equivalent:
\begin{itemize}
 \item[(a)] $F$ has finite projective dimension in $\CA$
 \item[(b)] $F$ has projective dimension at most one
 \item[(c)] $F$ has finite injective dimension
 \item[(d)] $F$ has injective dimension at most one
 \item[(e)] $F\in\CE$
 \item[(f)] $F$ is in the image of $h$.
\end{itemize}
Next, we say that a map $v$ in $\CC$ is \emph{phantom} if $h(v)=0$.
In a Brown category, the composite of any two phantom maps is zero, so
the phantoms form a square-zero ideal.  (Benson~\cite{be:pmpiii} has
shown that this can fail when $\CC$ is not a Brown category; in
particular, it fails when $\CC=\Stab_{kG}$, $k$ is an uncountable
field of characteristic $p$, and the $p$-rank of $G$ is at least two.)
Some further properties of phantom maps are studied
in~\cite{chho:pmc}.

Finally, consider a diagram $X\:I\xra{}\CC$, where $I$ is a filtered
category.  A \emph{weak colimit} for the diagram consists of a object
$U$ and compatible maps $X_i\xra{}U$ for $i\in I$, such that the
induced map $[U,Y]\xra{}\invlim_I[X_i,Y]$ is surjective for all
$Y\in\CC$.  Such a weak colimit is \emph{minimal} if the induced map
$\colim_i[Z,X_i]\xra{}[Z,U]$ is a bijection for all $Z\in\CF$.  In a
Brown category, it is known that
\begin{itemize}
 \item[(a)] Every filtered diagram of small objects has a minimal weak
  colimit, which is a retract of any other weak colimit.
 \item[(b)] Every object can be expressed as the minimal weak colimit
  of a filtered diagram of small objects.
\end{itemize}


\end{document}